\documentclass[11pt,oneside]{amsart} 
\usepackage{amsmath,graphicx}
\newtheorem{thm}{Theorem}[section]
\newtheorem{lem}[thm]{Lemma}
\newtheorem{cor}[thm]{Corollary}
\theoremstyle{definition}
\newtheorem{defn}[thm]{Definition}
\newtheorem{rem}[thm]{Remark}
\newtheorem{exmp}[thm]{Example}
\newtheorem{conv}[thm]{Convention}
\newtheorem{prob}[thm]{Open Question}
\newcommand{\Q}{\ensuremath{\mathbb{Q}}}
\newcommand{\R}{\ensuremath{\mathbb{R}}}
\newcommand{\euc}{\ensuremath{\mathbb{E}}}
\newcommand{\sph}{\ensuremath{\mathbb{S}}}
\newcommand{\hyp}{\ensuremath{\mathbb{H}}}
\newcommand{\script}[1]{\ensuremath{\mathcal{#1}}}
\newcommand{\C}{\script{C}}

\newcommand{\G}{\script{G}}
\newcommand{\eS}{\script{S}}
\newcommand{\smallcaps}[1]{\textrm{\textsc{#1}}}
\newcommand{\cat}{\smallcaps{CAT}}

\newcommand{\link}{\smallcaps{link}}
\newcommand{\shapes}{\smallcaps{Shapes}}
\newcommand{\length}{\textit{length}}

\begin{document}
%%%%%%%%%%%%%%%%

\title{CAT(0) is an algorithmic property}

\author[M.~Elder]{Murray Elder}
\address{Dept. of Mathematics\\
	Tufts University\\
	Medford, MA 02155}
\email{Murray.Elder@tufts.edu}

\author[J.~McCammond]{Jon McCammond$\ \!{ }^1$}
\address{Dept. of Math.\\
         University of California, Santa Barbara\\
         Santa Barbara, CA 93106}
\email{jon.mccammond@math.ucsb.edu}

\date{\today}

\begin{abstract}
In this article we give an explicit algorithm which will determine, in
a discrete and computable way, whether a finite piecewise Euclidean
complex is non-positively curved.  In particular, given such a complex
we show how to define a boolean combination of polynomial equations
and inequalities in real variables, i.e. a real semi-algebraic set,
which is empty if and only if the complex is non-positively curved.
Once this equivalence has been shown, the main result follows from a
standard theorem in real algebraic geometry.
\end{abstract}

\footnotetext[1]{Supported under NSF grant no.\ DMS-99781628}
\maketitle

\tableofcontents

In this article we prove that for every finite collection of Euclidean
polytopes, there exists a finite list of forbidden configurations
which determine the non-positively curved complexes that can be built
out of these shapes.  More specifically, a complex built out of these
shapes will be non-positively curved if and only if it avoids these
forbidden configurations.  The general result for arbitrary curvatures
is as follows:

\renewcommand{\thethm}{\Alph{thm}}
\begin{thm}\label{thm:minors}
If $\eS$ is a finite collection of shapes with curvature $\kappa$,
then there exists a finite list of configurations $\C$ such that an
$M_\kappa$-complex $K$ with $\shapes(K) \subset \eS$ is locally
$\cat(\kappa)$ if and only if $K$ avoids all of the configurations in
$\C$.
\end{thm}
\renewcommand{\thethm}{\thesection.\arabic{thm}}

Notice, however, that Theorem~\ref{thm:minors} merely asserts the
existence of a finite list of forbidden configurations that need to be
avoided.  It does not, in and of itself, show that there is a
procedure for deciding whether a particular finite complex is locally
$\cat(\kappa)$.  The second half of the article is devoted to proving
this stronger assertion.  In particular we prove the following:

\renewcommand{\thethm}{\Alph{thm}}
\begin{thm}\label{thm:main}
There exists an algorithm which determines whether or not a finite
$M_\kappa$-complex is locally $\cat(\kappa)$.
\end{thm}
\renewcommand{\thethm}{\thesection.\arabic{thm}}

Two special cases are of particular note.  Theorem~\ref{thm:main}
shows that it is possible to test whether a finite, piecewise
Euclidean complex is non-positively curved and to test whether a
finite, piecewise hyperbolic complex is negatively curved.

The proof of Theorem~\ref{thm:main} proceeds by converting curvature
considerations into a real semi-algebraic set which is empty if and
only if the complex is $\cat(\kappa)$.  A \emph{real semi-algebraic
subset of $\R^n$} is a subset which can be described by a finite
boolean combination of polynomial equations and inequalities in $n$
real variables.  The main property of real semi-algebraic sets we will
need is that given a boolean combination of polynomial equations and
inequalities in $n$ variables there is an algorithm which decides
whether this system has a solution.  This is a special case of
Tarski's theorem that the elementary theory of the reals is decidable.
See \cite{BoCoRo} or \cite{CaJo} for details.  Thus, once the
conversion to real algebraic geometry has been completed, the
existence of an algorithm will be immediate.

Finally, we should note that in dimension~$3$ there exists a
particularly elementary procedure which avoids real algebraic geometry
altogether.  Since this algorithm differs substantially in flavor from
the one given below, it will be described elsewhere.  See
\cite{ElMc-three} for details.

\vspace{1em}
%%%%%%%%%%%%%%%%%%%%%%%%%%%%%%%%%%%%%%%%%%
\noindent\textbf{Overview of the Sections}\\
%%%%%%%%%%%%%%%%%%%%%%%%%%%%%%%%%%%%%%%%%%
In the first two sections we briefly review some basic definitions and
results from polyhedral geometry, and comparison geometry.  They are
included for completeness and can readily be skipped by the reader
familiar with these areas.

In Section~\ref{sec:galleries} we introduce our main technical tool,
that of a gallery, and establish its main properties.  In
Section~\ref{sec:piecewise} we analyze how to describe piecewise
geodesics in galleries using polynomial equations and inequalities and
in Section~\ref{sec:geos} we perform a similar analysis for local
geodesics.  These two sections represent the heart of the proof.

In Section~\ref{sec:coefficients} we show that the coefficients that
arise in the encoding process are no worse that the original lengths
in the original complex, and finally, in Section~\ref{sec:main} we
present the proof of the main theorem, Theorem~\ref{thm:main}.

%%%%%%%%%%%%%%%%%%%%%%%%%%%%%%%%%%%%%%%%%%%%%%%%%%%
\section{Polyhedral geometry}\label{sec:polyhedral}
%%%%%%%%%%%%%%%%%%%%%%%%%%%%%%%%%%%%%%%%%%%%%%%%%%%
We begin by reviewing the necessary background regarding polyhedral
geometry.  A metric space $(K,d)$ is called a \emph{geodesic metric
space} if every pair of points in $K$ can be connected by a geodesic.
The geodesic metric spaces we will be primarily interested in will be
cell complexes constructed out of convex polyhedral cells in $\hyp^n$,
$\euc^n$ or $\sph^n$.  After describing these spaces, we will quickly
review their relationship with comparison geometry.  Details can be
found in \cite{Ba}, \cite{Br-thesis}, \cite{BrHa}, and \cite{Mo}.

\begin{defn}[Polyhedral cells]\label{def:cells}
A \emph{convex polyhedral cell} in $\hyp^n$ or $\euc^n$is the convex
hull of a finite set of points.  The convex hull of $n+1$ points in
general position is an \emph{$n$-simplex}. A \emph{polyhedral cone} in
$\euc^n$ is the positive cone spanned by a finite set of vectors.  If
the original vectors are linearly independent, it is a
\emph{simplicial cone}.  A \emph{cell (simplex) in $\sph^n$} is the
intersection of a polyhedral cone (simplicial cone) in $\euc^n$ with
$\sph^n$.  

A spherical cell which does not contain a pair of antipodal points is
\emph{proper}.  All Euclidean and hyperbolic cells are considered
proper.  Notice that every spherical cell can be subdivided into
proper spherical cells by it cutting along the coordinate axes.  If
$\sigma$ denotes a proper convex polyhedral cell, then $\sigma^\circ$
will denote its interior and $\partial \sigma$ will denote its
boundary.
\end{defn}

\begin{defn}[$M_\kappa$-complexes]\label{def:m-complexes}
An \emph{$\hyp$-complex} [$\euc$-complex, $\sph$-complex] is a
connected cell complex $K$ make up of proper hyperbolic [Euclidean,
spherical] cells glued together by isometries along faces.  A cell
complex which has an $\hyp$-complex, an $\euc$-complex, or an
$\sph$-complex structure, will be called a \emph{metric polyhedral
complex}, or \emph{$M_\kappa$-complex} for short where $\kappa$
denotes the curvature constant common to all of its cells.  More
generally, an $M_\kappa$-complex is formed from polyhedral cells with
constant curvature $\kappa$.
\end{defn}

\begin{conv}[Subdivisions]\label{conv:subdivisions}
As noted above, every spherical cell can be subdivided into proper
spherical cells.  Since proper spherical cells are required in an
$\sph$-complex, a subdivision will occasionally be necessary in order
to convert a complex built out of pieces of spheres to be considered
an $\sph$-complex.
\end{conv}

\begin{defn}[Shapes]\label{def:shapes}
If $K$ is an $M_\kappa$-complex, then the isometry types of the cells
of $K$ will be called the \emph{shapes of $K$} and the collection of
these isometry types will be denoted $\shapes(K)$.  When $\shapes(K)$
is finite, $K$ is said to have only \emph{finitely many shapes}.
Notice that since cells of different dimensions necessarily have
different isometry types, finitely many shapes implies that $K$ is
finite dimensional.  It does not, however, imply that $K$ is locally
finite.  Notice also that if $K$ is an $M_\kappa$-complex with only
finitely many shapes, then there is a subdivision $K'$ of $K$ where
the cells of $K'$ are proper and simplicial, and $\shapes(K')$ remains
finite.
\end{defn}

\begin{defn}[Paths and Loops]
A \emph{path $\gamma$} in a metric space $K$ is a continuous map
$\gamma:[0,\ell]\rightarrow K$.  A path is \emph{closed} if $\gamma(0)
= \gamma(\ell)$.  A \emph{loop} is a closed path where the basepoint
$\gamma(0)$ has been forgotten.  Technically, a loop is viewed as a
continuous map from a circle to $K$.  The circle of perimeter $\ell$
will be denoted $C_\ell$.
\end{defn}

\begin{defn}[Piecewise geodesics]\label{def:metric}
Let $K$ be an $M_\kappa$-complex. A \emph{piecewise geodesic} $\gamma$
in $K$ is a path $\gamma:[a,b] \rightarrow K$ where $[a,b]$ can be
subdivided into a finite number of subintervals so that the
restriction of $\gamma$ to each closed subinterval is a path lying
entirely in some closed cell $\sigma$ of $K$ and that this path is the
unique geodesic connecting its endpoints in the metric of $\sigma$.
The \emph{length} of $\gamma$, denoted $\length(\gamma)$, is the sum
of the lengths of the geodesics into which it can be partitioned.  A
\emph{closed piecewise geodesic} and a \emph{piecewise geodesic loop}
are defined similarly.  The \emph{intrinsic metric} on $K$ is defined
as follows:
\[ d(x,y) = \inf \{\length(\gamma) | \textrm{$\gamma$ is a piecewise
geodesic from $x$ to $y$} \}\]
In general $d$ is only a pseudometric, but when $K$ has only finitely
many shapes, $d$ is a well-defined metric and $(K,d)$ is a geodesic
metric space \cite[Theorem~7.19]{BrHa}.  
\end{defn}

\begin{conv}[Parameterizations]
If $\gamma:[0,\ell]$ is a piecewise geodesic, we will always assume
that the map $\gamma$ has been reparameterized by arc-length.  In
particular, $\ell$ should be the length of $\gamma$ and for all
subintervals $[a,b]$ in $[0,\ell]$, the length of $\gamma([a,b])$
should be exactly $b-a$.  Following the same convention, if a loop is
a closed piecewise geodesic of length $\ell$, then we will assume that
its domain is a circle whose circumference is $\ell$.  We will use
$C_\ell$ to denote such a circle and we will identify the points on
$C_\ell$ with reals mod $\ell$ and in particular with the points
$[0,\ell)$.
\end{conv}

\begin{defn}[Size]
The \emph{size} of a piecewise geodesic $\gamma:[0,\ell]\rightarrow K$
is the number of open cells of $K$ through which $\gamma([0,\ell])$
passes, with multiplicities.  Technically, the size of $\gamma$ is the
minimal number of subintervals (open, half-open, or closed) into which
$[0,\ell]$ must partitioned so that the image of each subinterval lies
in a single open cell of $K$.  Note that some of these subintervals
may be single points.  The fact that $\gamma$ is a piecewise geodesic
ensures that the size of $\gamma$ is finite.
\end{defn}

\begin{defn}[Links]\label{def:links}
Let $K$ be an $M_\kappa$-complex with only finitely many
shapes and let $x$ be a point in $K$.  The set of unit tangent vectors
to $K$ at $x$ is naturally an $\sph$-complex called the \emph{link of
$x$ in $K$}, or $\link(x,K)$.  If $K$ has only finitely many shapes,
then $\link(x,K)$ has only finitely many shapes.

When $x$ lies in the interior of a cell $B$ of $K$, $\link(x,B)$ is a
sphere of dimension $k = \dim B - 1$ sitting inside $\link(x,K)$.
Moreover, the complex $\link(x,K)$ can be viewed as a spherical join
of $\sph^k$ and another $\sph$-complex, denoted $\link(B,K)$, which
can be thought of as the unit tangent vectors to $x$ in $K$ which are
orthogonal to $B$.  The complex $\link(B,K)$ is called the \emph{link
of the cell $B$ in $K$}.  Once again, if $K$ has only finitely many
shapes, then $\link(B,K)$ has only finitely many shapes as well.
\end{defn}

\begin{defn}[Local geodesics]
Let $K$ be an $M_\kappa$-complex. A piecewise geodesic $\gamma$ in $K$
is called a \emph{local geodesic} if for each point $x$ on $\gamma$,
the incoming and outgoing unit tangent vectors to $\gamma$ at $x$ are
at a distance of at least $\pi$ from each other in $\link(x,K)$.
\end{defn}

The size and length of local geodesics are closely related.

\begin{thm}[Bridson]\label{thm:quasi-geo}
If $K$ is an $M_\kappa$-complex with only finitely many shapes, then
for every $\ell > 0$ there exist an integer $N > 0$, depending only on
$\shapes(K)$, such that every local geodesic of size at least $N$ has
length at least $\ell$.
\end{thm}

For a proof of this result see \cite[Theorem~1.11]{Br-thesis} or
\cite[Theorem~I.7.28]{BrHa}.  Both references state this theorem in a
more general form using the concept of a taut-$m$-string.  We merely
need to note that local geodesics are taut-$m$-strings and the size of
a local geodesic $\gamma$, as defined above, is within a constant
factor of the integer $m$ when $\gamma$ is viewed as a
taut-$m$-string.

\begin{rem}[Constructing $N$]\label{rem:N}
The \emph{existence} of a constant such as $N$ is important, but for
our purposes, we need to know more.  There needs to be an algorithm
which constructs an integer $N$ which will work solely from
$\shapes(K)$ and the real number $\ell$.  There is indeed such an
algorithm, and we will return to this issue later in the article.  See
Lemma~\ref{lem:N}.

Notice that it is in fact sufficient to construct an $N$ that will
work for one particular value of $\ell$.  That is, if every local
geodesic of size $N$ has length at least $\ell$, then every local
geodesic of size $k\cdot N$ will have length at least $k\cdot \ell$.
Thus for an arbitrary $\beta$, the least integer greater than
$\frac{\beta}{\ell}\cdot N$ will work for $\beta$.
\end{rem}

%%%%%%%%%%%%%%%%%%%%%%%%%%%%%%%%%%%%%%%%%%%%%%%%%%%
\section{Comparison geometry}\label{sec:comparison}
%%%%%%%%%%%%%%%%%%%%%%%%%%%%%%%%%%%%%%%%%%%%%%%%%%%
Metric polyhedral complexes are particularly useful in the creation of
metric spaces of non-positive curvature.  As in the previous section,
details can be found in \cite{Ba}, \cite{Br-thesis}, \cite{BrHa}, and
\cite{Mo}.

\begin{defn}[Globally $\cat(\kappa)$]
Let $K$ be a geodesic metric space, let $T$ be a geodesic triangle in
$K$, and let $\kappa = -1$, [or $0$, or $1$].  A \emph{comparison
triangle} for $T$ is a triangle $T'$ in $\hyp^2$, [or $\euc^2$ or
$\sph^2$] with the same side lengths as $T$.  Notice that for every
point $x$ on $T$, there is a corresponding point $x'$ on $T'$.  The
space $K$ is called \emph{globally $\cat(\kappa)$}, if for any
geodesic triangle $T$ in $K$ [of perimeter less than $2\pi$ when
$\kappa=1$] and for any points $x$ and $y$ on $T$, the distance from
$x$ to $y$ in $K$ is less than or equal to the distance from $x'$ to
$y'$ in $\hyp^2$ [or $\euc^2$ or $\sph^2$].  Finally, a space $K$ is
called \emph{locally $\cat(\kappa)$} if every point in $K$ has a
neighborhood which is globally $\cat(\kappa)$.  Locally $\cat(0)$
spaces are often referred to as \emph{non-positively curved} and
locally $\cat(-1)$ spaces are called \emph{negatively curved}.
\end{defn}

\begin{thm}
In a globally $\cat(0)$ space, every pair of points is connected by a
unique geodesic, and a path is a geodesic if and only if it is a local
geodesic.
\end{thm}

The next two results about $M_\kappa$-complexes show how global
properties such as $\cat(\kappa)$ can be reduced to local properties
and how local properties can be reduced to the existence of geodesics
in $\sph$-complexes.

\begin{thm}\label{thm:local-to-global}
Let $K$ be an $M_\kappa$-complex which contains only finitely many
shapes.
\begin{enumerate}
\item If $K$ is an $\hyp$-complex, then $K$ is globally $\cat(-1)$ if
and only if it is locally $\cat(-1)$ and simply-connected.
\item If $K$ is an $\euc$-complex, then $K$ is globally $\cat(0)$ if
and only if it is locally $\cat(0)$ and simply-connected.
\item If $K$ is an $\sph$-complex, then $K$ is globally $\cat(1)$ if
and only if it is locally $\cat(1)$ and and there are no geodesic
cycles of length $< 2\pi$.
\end{enumerate}
\end{thm}

\begin{thm}\label{thm:local-via-geos}
If $K$ is an $M_\kappa$-complex, then the following are equivalent:
\begin{enumerate}
\item $K$ is locally $\cat(\kappa)$.
\item The link of each vertex in $K$ is globally $\cat(1)$.
\item The link of each cell of $K$ is an $\sph$-complex which
contains no closed geodesic cycle of length less than $2\pi$.
\end{enumerate}
\end{thm}

Thus showing that $\hyp$-complexes are $\cat(-1)$ or that
$\euc$-complexes are $\cat(0)$ ultimately depends on being able to
show that various $\sph$-complexes contain no short geodesic cycles.
In order to study short geodesic cycles we introduce the concept of a
configuration.  Configurations are related to the ``finite models''
used by Bridson and Haefliger \cite{BrHa}.  We will use configurations
and the following theorem, which is a restatement of
Theorem~\ref{thm:quasi-geo}, to prove Theorem~\ref{thm:minors}.

\begin{thm}\label{thm:finite-subcomplex}
If $K$ is an $\sph$-complex with only finitely many shapes, then there
is a constant $N$ depending only on $\shapes(K)$, such that every
local geodesic of length less than $2\pi$ has size less than $N$.  In
other words, if $K$ contains a short closed geodesic $\gamma$, then
$\gamma$ is contained in a finite subcomplex of $K$ which is the union
of at most $N$ cells from $\shapes(K)$.
\end{thm}

\begin{defn}[Configurations] 
A \emph{configuration} $C$ is a finite $\sph$-complex $C$ which
contains at least one closed geodesic $\gamma$ of length less than
$2\pi$.  An $M$-complex $K$ \emph{contains a configuration} $C$ if $C$
isometrically embeds as a subcomplex of $\link(B,K)$ for some cell $B$
in $K$ and under this embedding at least one of the short closed
geodesics contained in $C$ is sent to a short closed geodesic in
$\link(B,K)$.  If $K$ does not contain a configuration $C$, then $K$
\emph{avoids} this configuration.
\end{defn}

If we fix the set of shapes under consideration and let the
$M_\kappa$-complex vary, then Theorem~\ref{thm:finite-subcomplex} can
be restated as an assertion about the existence of finite lists of
forbidden configurations which characterize which of the complexes
built out of these shapes are $\cat(\kappa)$.

\renewcommand{\thethm}{\ref{thm:minors}}
\begin{thm}
If $\eS$ is a finite collection of shapes with curvature $\kappa$,
then there exists a finite list of configurations $\C$ such that an
$M_\kappa$-complex $K$ with $\shapes(K) \subset \eS$ is locally
$\cat(\kappa)$ if and only if $K$ avoids all of the configurations in
$\C$.
\end{thm}
\renewcommand{\thethm}{\thesection.\arabic{thm}}
\addtocounter{thm}{-1}

\begin{proof}
By Theorem~\ref{thm:local-via-geos} it is sufficient to determine
whether the link of a cell in $K$ contains a closed geodesic of length
less than $2\pi$, and by Theorem~\ref{thm:finite-subcomplex} this
closed geodesic must live in a finite subcomplex whose size is
uniformly bounded by a constant which depends only on the set of
shapes in $\eS$.  Since there are only a finite number of finite
complexes of bounded size built out of the links of faces of shapes in
$\eS$ (up to isometry), the list of configurations which will destroy
the property of being locally $\cat(\kappa)$ is contained in a finite
list of possibilities.  Even though this does not determine which ones
they are, this does show that the list of forbidden configurations is
finite.
\end{proof}

As we noted in the introduction, Theorem~\ref{thm:minors} merely
asserts the existence of a finite list of forbidden configurations
that need to be avoided but it does not provide an algorithm to
construct the list.  At present there are very few instances where the
finite list for some collection of shapes is known explicitly.  One of
the main results along these lines is Moussong's Lemma which was first
stated by Gromov (\cite{Gr}) in the special case where all $1$-cells
in $K$ have length $\pi/2$ and proved in general by Moussong in his
dissertation \cite{Mo}.

\begin{exmp}[Moussong's Lemma]
Moussong's Lemma states that if $K$ is a simplicial $\sph$-complex
with only finitely many shapes and all of the $1$-cells of $K$ have
length at least $\frac{\pi}{2}$, then $K$ is globally $\cat(1)$ if and
only if $K$ is a metric flag complex.  A simplicial complex is called
a \emph{flag complex} if every $1$-skeleton of an $n$-simplex is
filled in with an $n$-simplex.  A simplicial $\sph$-complex is called
a \emph{metric flag complex} if every $1$-skeleton of an $n$-simplex
whose edge lengths are those of a possible proper spherical
$n$-simplex is filled in with a copy of that $n$-simplex.  This means,
for example, that  three $1$-cells in an $\sph$-complex each with
length $\pi$ need not and cannot bound a spherical triangle.

When all of the $1$-cells have length $\pi/2$, the metric flag
condition reduces to the (non-metric) flag condition.  Moussong's
lemma thus enables a researcher to easily determine whether an
$n$-dimensional cubical complex with the usual metric has a locally
$\cat(0)$ structure.  The finite list of forbidden configurations in
this case is the collection of empty simplices in dimensions $2$ up to
$n-1$.
\end{exmp}

One way to paraphrase Theorem~\ref{thm:minors} is that for any finite
collection of piecewise Euclidean or piecewise hyperbolic shapes,
there is a finite description (similar to that in Moussong's lemma) of
the combinatorial conditions needed in order for a complex built out
of these shapes to be $\cat(0)$ or $\cat(-1)$.  Such finite
descriptions must exist by Theorem~\ref{thm:minors}, but most of them
remain undiscovered.

%%%%%%%%%%%%%%%%%%%%%%%%%%%%%%%%%%%%%%%%
\section{Galleries}\label{sec:galleries}
%%%%%%%%%%%%%%%%%%%%%%%%%%%%%%%%%%%%%%%%
In this section we define one of our central concepts: that of a
gallery.  A gallery is a more precise tool than a configuration, as
will become clear below.  The term ``gallery'' has been borrowed from
the study of Coxeter groups since a geodesic gallery in a Coxeter
complex will be an example of a linear gallery as defined below.  The
galleries defined here will be either linear or circular and we begin
by defining the linear ones.

\begin{defn}[Linear Galleries]\label{def:lin-gallery}
Let $K$ be a $M_\kappa$-complex, let $\gamma:[0,\ell]\rightarrow K$ be
a piecewise geodesic in $K$, let $k$ is the size of $\gamma$, and let
$\{\sigma_i^\circ\}_{i=1}^{k}$ be the sequence of open cells that
$\gamma$ passes through.  Finally, assume that for all $1<i<k$,
$\sigma_i$ is either a common face of both $\sigma_{i+1}$ and
$\sigma_{i-1}$ which are distinct or else both of these are
incompatible faces of $\sigma_i$.  In the former case, we say that
$\sigma_i$ is a \emph{bottom cell}, and in the latter $\sigma_i$ is an
\emph{top cell} in the sequence.

The \emph{linear gallery $\G$ determined by $\gamma$} is constructed
inductively from this sequence of open cells.  To start the induction,
let $\G_1$ be a copy of the closed cell $\sigma_1$, let $\alpha_1$ be
all of $\G_1$, and let $\phi_1: \alpha_1 \rightarrow \sigma_1$ be a
specific isometry between $\alpha_1$ and $\sigma_1$.  The cell
$\alpha_i$ in $\G_i$, will always be isometric to $\sigma_i$ by a
specific isometry $\phi_i$ and the complex $\G_{i+1}$ will always
contain the complex $\G_i$ as a subcomplex.  We call $\alpha_i$ the
\emph{active cell in $\G_i$ }.  Next, assume that for some $i \geq 1$
we have defined the complex $\G_i$, its active cell $\alpha_i$, and an
isometric $\phi_i:\alpha_i \rightarrow \sigma_i$.  The inductive step
depends on whether $\sigma_i$ is a top cell or a bottom cell.
\begin{itemize}
\item
If $\sigma_{i+1}$ is a bottom cell, then define $\G_{i+1} = \G_i$,
define $\alpha_{i+1} = \phi_i^{-1}(\sigma_{i+1})$ and define $\phi_{i+1}$
as the restriction of $\phi_i$ to $\alpha_{i+1}$.
\item
If $\sigma_i$ is a top cell, then define let $\alpha_{i+1}$ be a copy
of $\sigma_{i+1}$, let $\phi_{i+1}:\alpha_{i+1} \rightarrow
\sigma_{i+1}$ be an isometry between the two, and let $\G_{i+1}$ be
the complex formed by gluing $\G_i$ to $\alpha_{i+1}$ along their
distinguished faces isometric to $\sigma_{i}$.  Specifically identify
$\phi_{i+1}^{-1}(\sigma_i) \subset \alpha_{i+1}$ with $\alpha_i
\subset \G_i$ in the obvious fashion.  Note that when $\G_i$ is viewed
as a subcomplex of $\G_{i+1}$, $\phi_{i+1}$ is an extension of
$\phi_i$.
\end{itemize}

The linear gallery $\G$ is simply the final complex $\G_k$.  The cells
$\alpha_1$ and $\alpha_k$ in $\G$ are the \emph{endcells} of $\G$.
The \emph{interior of $\G$}, denoted $\G^\circ$, is the union of the
\emph{open} cells $\alpha_i$, $i=1,\ldots, k$ in $\G$.  The
\emph{boundary of $\G$}, denoted $\partial \G$, is the open cells in
$\G$ which are not in the interior.  Notice that since the various
maps $\phi_i$ are compatible on their overlaps, the gallery $\G$ comes
equipped with a map $\phi:\G \rightarrow K$ which is an isometry when
restricted to any of the closed cells of $\G$.
\end{defn}

\begin{figure}
\begin{tabular}{cc}
\begin{tabular}{c}\includegraphics{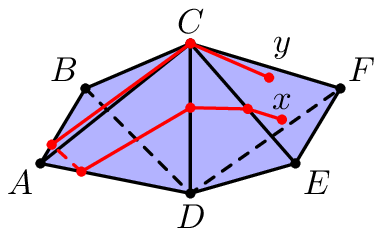}\end{tabular}& 
\begin{tabular}{c}\includegraphics{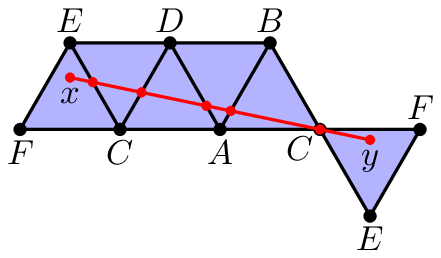}\end{tabular}\\
\begin{tabular}{c}\includegraphics{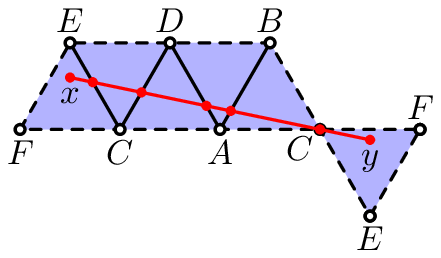}\end{tabular}& 
\begin{tabular}{c}\includegraphics{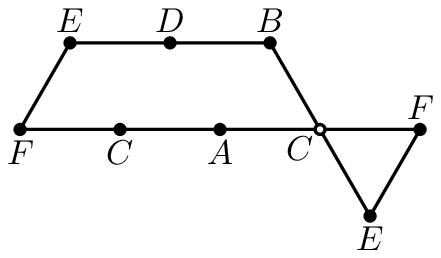}\end{tabular}
\end{tabular}
\caption{The $2$-complex and linear gallery described in
Example~\ref{exmp:2-gallery}.
\label{fig:2-gallery}}
\end{figure}

\begin{exmp}\label{exmp:2-gallery}
Let $K$ be the $2$-dimensional $\euc$-complex formed by attaching the
boundaries of two regular Euclidean tetrahedra along a $1$-cell.  The
complex $K$ is shown in the upper left corner of
Figure~\ref{fig:2-gallery}.  Let $\gamma$ be the geodesic shown which
starts at $x$ travels across the front of $K$, around the back, over
the top, and ends at $y$.  The gallery determined by $\gamma$ in shown
in the upper right corner, its interior in the lower left corner and
its boundary in the lower right of Figure~\ref{fig:2-gallery}.
\end{exmp}

\begin{defn}[Circular Galleries]\label{def:cir-gallery}
Let $K$ be a $M_\kappa$-complex, let $\gamma:C_\ell \rightarrow K$ be
a closed piecewise geodesic in $K$, let $k$ is the size of $\gamma$,
and let $\{\sigma_i^\circ\}_{i=1}^{k}$ be the sequence of open cells
that $\gamma$ passes through.  Finally, assume that for all $1 \leq i
\leq k$, $\sigma_i$ is either a common face of both $\sigma_{i+1}$ and
$\sigma_{i-1}$ which are distinct or else both of these are
incompatible faces of $\sigma_i$ where the subscripts are interpreted
$\mod k$..  In the former case, we say that $\sigma_i$ is a
\emph{bottom cell}, and in the latter $\sigma_i$ is an \emph{top cell}
in the sequence.

As in Definition~\ref{def:lin-gallery}, the \emph{circular gallery
$\G$ determined by $\gamma$} is constructed inductively from the
sequence of open cells that $\gamma$ passes through.  Let $a$ be an
arbitrary point in $C_\ell$, say one whose image lies in
$\sigma_1^\circ$.  If we cut the circle at $a$ we obtain a path
$\gamma'$ of length $\ell$ in $K$ which is a piecewise geodesic and
which happens to start and end at the same point, $\gamma'(a)$.  Next,
we construct a linear gallery $\G'$ for the path $\gamma'$ as above.
Finally, we identify the endcells of $\G'$ to form the circular
gallery $\G$.

It should be clear that the end result is independent of the choice of
$a$, and that the circular gallery $\G$ comes equipped with a map
$\phi:\G\rightarrow K$ which is an isometry when restricted to any of
the cells of $\G$.
\end{defn}

\begin{defn}[Galleries in complexes]
Let $K$ be a $M_\kappa$-complex with only finitely many shapes.  A
complex $\G$ with a map $\phi:\G\rightarrow K$ will be called a
\emph{linear gallery in $K$} if there is a piecewise geodesic $\gamma$
satisfying the necessary restrictions which determines $\G$ and
$\phi$.  A \emph{circular gallery in $K$} is defined analogously.
We should note that when we speak of a linear gallery $\G$ in $K$, we
implicitly assume that $\G$ has a specified pair of endcells, $\sigma$
and $\tau$, and a specified map $\phi:\G \rightarrow K$.  From this
information we can recover the interior $\G^\circ$ and the linear
ordering on the open cells in $\G^\circ$.  Similarly, when we speak of
a circular gallery $\G$ in $K$, we implicitly assume that $\G$ has a
specified map $\phi:\G \rightarrow K$.  From $\G$ we can reconstruct
its interior $\G^\circ$ and the circular ordering of open cells in
$\G^\circ$.
\end{defn}

\begin{rem}[Recognizing Galleries]\label{rem:recognize}
Let $\phi:\G\rightarrow K$ be map from a finite complex to $K$ which
is an isometry when restricted to a closed cell of $\G$, and let
$\sigma$ and $\tau$ be two distinguished closed cells in $\G$.  It is
easy to determine whether $\G$ is a linear gallery with endcells
$\sigma$ and $\tau$ as determined by some piecewise geodesic in $K$.
If $\G$ is a gallery, then the maximal cells in the face lattice of
$\G$ are the top cells of the gallery.  There needs to exist an
ordering of these cells so that adjacent top cells have non-trivial
intersections. These intersections are the bottom cells. If $\sigma$
or $\tau$ is not in the list of top cells then the cells containing
them are also bottom cells. At this point, we simply check whether the
two bottom cells adjacent to a top cell are incompatible.  If all of
these conditions are true, then $\G$ is a gallery determined by some
piecewise geodesic.  Since there are only finitely many maximal cells
in $\G$ and only finitely many linear ordering on this set, we can
check all of the possibilities.  A similar procedure will determine
whether a map $\phi:\G\rightarrow K$ from a finite complex to $K$
which is an isometry when restricted to a closed cell of $\G$ is a
circular gallery determined by some closed piecewise geodesic.
\end{rem}

\begin{lem}\label{lem:lifts}
Let $K$ be an $M_\kappa$-complex and let $\phi:\G\rightarrow K$ be the
linear gallery determined by a piecewise geodesic $\gamma:[0,\ell]
\rightarrow K$.  The map $\gamma$ lifts through $\phi$ in the sense
that there is path $\gamma':[0,\ell]\rightarrow \G$ with $\gamma =
\phi \circ \gamma'$.  Moreover, $\gamma'$ is a piecewise geodesic in
$\G$ which starts in the one endcell of $\G$, ends in the other, and
passes monotonically through the open cells in $\G^\circ$.
\end{lem}

\begin{proof}
The notation of Definition~\ref{def:lin-gallery} will be used without
comment.  Partition the path $\gamma$ into open, closed and half-open
subintervals so that the $i$-th subinterval is sent to the $i$-th open
cell, $\sigma_i^\circ$, in $K$ that $\gamma$ passes through.  If
$\sigma_i$ is a top cell, then it contains the adjacent bottom cells
in its boundary and thus $\sigma_i$ contains the images of the
$(i-1)$-st and the $(i+1)$-st subintervals if they exist.  Thus all
three subintervals can be simultaneously lifted to $\alpha_i \cap
\G^\circ$ via the isometry $\phi_i^{-1}$.  Since these lifts agree on
the bottom cells where they overlap, they combine to form the required
path $\gamma'$.  The second assertion is immediate from the definition
of $\G$.
\end{proof}

\begin{rem}[Geodesics in galleries]\label{rem:geos-in-galleries}
Note that several piecewise geodesics may determine the same linear
gallery $\G$ in $K$.  If $\gamma$ is one of these paths which
determine $\G$, then we say that $\G$ \emph{contains} $\gamma$.
Notice in particular that if $\G$ contains $\gamma$, then the lift
$\gamma'$ of $\gamma$ to $\G$ will be a piecewise geodesic in
$\G^\circ$ which starts in one endcell of $\G$, ends in the other
endcell of $\G$ and which progresses monotonically from one endcell to
the other.  Conversely given a piecewise geodesic $\gamma'$ in $\G$
satisfying these restrictions, it is easy to see that its $\phi$-image
in $K$ will be a piecewise geodesic $\gamma$ which determines $\G$.
\end{rem}

\begin{lem}\label{lem:subintervals}
Let $K$ be an $M_\kappa$-complex, and let $\phi:\G\rightarrow K$ be
either the linear gallery determined by a geodesic $\gamma:[0,\ell]
\rightarrow K$, or the circular gallery determined by a closed
geodesic $\gamma:C_\ell \rightarrow K$.  If $\gamma'$ is the lift of
$\gamma$ to $\G$ and $\alpha_i^\circ$ is an open cell in $\G^\circ$,
then the inverse image of $\alpha_i^\circ$ under $\gamma'$ is a closed
subinterval of $[0,\ell]$ or $C_\ell$ when $\alpha_i$ is a bottom cell
and an open subinterval when $\alpha_i$ is a top cell.
\end{lem}

\begin{proof}
Notice that $\alpha_i^\circ$ is a closed subspace of $\G^\circ$ when
$\alpha_i$ is a bottom cell in $\G$ and an open subspace when
$\alpha_i$ is a top cell in $\G$.  Thus the preimage of
$\alpha_i^\circ$ is a closed [respectively open] subspace of
$[0,\ell]$ or $C_\ell$ when $\alpha_i$ is a bottom [top] cell.  Since
$\gamma'$ proceeds monotonically through the interiors of the
$\alpha_i$ (Lemma~\ref{lem:lifts}), these open/closed subspaces can
only be open/closed subintervals.
\end{proof}

\begin{lem}\label{lem:bottom-cells}
Let $K$ be an $M_\kappa$-complex, and let $\phi:\G\rightarrow K$ be
either the linear gallery determined by a geodesic $\gamma:[0,\ell]
\rightarrow K$, or the circular gallery determined by a closed
geodesic $\gamma:C_\ell \rightarrow K$.  If $\gamma'$ is the lift of
$\gamma$ to $\G$ and $\alpha$ is a bottom cell in $\G$, then the
inverse image of $\alpha^\circ$ under $\gamma'$ is a single point.
The same is true for circular gallery.
\end{lem}

\begin{proof}
If $\alpha$ is a bottom cell whose preimage under $\gamma'$ is not a
single point, then by considering $\alpha$ and an adjacent top cell,
it is easy to see that $\gamma'$ (and hence $\gamma$) is not a local
geodesic.  This shows that preimages of the interiors of bottom cells
must be single points.
\end{proof}

\begin{lem}\label{lem:lin-determined}
Let $K$ be an $M_\kappa$-complex and let $\phi:\G\rightarrow K$ be the
linear gallery determined by a geodesic $\gamma:[0,\ell] \rightarrow
K$.  If $\gamma'$ is the lift of $\gamma$ to $\G$, and $L = \{0=x_0 <
x_1 < \cdots < x_k=\ell\} \subset [0,\ell]$ is the list of point
preimages of bottom cells under $\gamma'$ together with the endpoints
of the the interval, then the image of $\gamma'$ in $\G$ is completely
determined by the images of the points in $L$.  Moreover, $\gamma'$ is
an embedding.
\end{lem}

\begin{proof}
If $\alpha$ denotes a top cell in $\G$, then the preimage of $\alpha$
under $\gamma'$ is an closed subinterval $[x_i,x_{i+1}]$ for some $i$
(Lemma~\ref{lem:subintervals} and Lemma~\ref{lem:bottom-cells}).
Since $\gamma'$ is a local geodesic, and since the image of this
subinterval lies completely in the closed cell $\alpha$, this image
must be the unique geodesic in $\alpha$ connecting $\gamma'(x_i)$ to
$\gamma'(x_{i+1})$.  This shows that $\gamma'$ is an embedding and
that it image is completely determined by the image of $L$ in $\G$ as
required.
\end{proof}

\begin{lem}\label{lem:cir-determined}
Let $K$ be an $M_\kappa$-complex and let $\phi:\G\rightarrow K$ be the
circular gallery determined by a closed geodesic $\gamma:C_\ell
\rightarrow K$.  If $\gamma'$ is the lift of $\gamma$ to $\G$, and $L
= \{0 \leq x_0 < x_1 < \cdots < x_k <\ell\} \subset [0,\ell)$ is the
list of point preimages of bottom cells under $\gamma'$, then the
image of $\gamma'$ in $\G$ is completely determined by the images of
the points in $L$.  Moreover, $\gamma'$ is an embedding.
\end{lem}

\begin{proof}
If $\alpha$ denotes a top cell in $\G$, then the preimage of $\alpha$
under $\gamma'$ is an closed subinterval $[x_i,x_{i+1}]$ for some $i$,
or $[x_k,x_0+\ell]$ (Lemma~\ref{lem:subintervals} and
Lemma~\ref{lem:bottom-cells}).  Since $\gamma'$ is a local geodesic,
and since the image of this subinterval lies completely in the closed
cell $\alpha$, this image must be the unique geodesic in $\alpha$
connecting $\gamma'(x_i)$ to $\gamma'(x_{i+1})$.  This shows that
$\gamma'$ is an embedding and that it image is completely determined
by the image of $L$ in $\G$ as required.
\end{proof}

\begin{lem}\label{lem:immersion}
Let $K$ be an $M_\kappa$-complex, let $\phi:\G\rightarrow K$ be the
linear gallery determined by a piecewise geodesic $\gamma:[0,\ell]
\rightarrow K$, and let $\gamma':[0,\ell] \rightarrow \G$ be the lift
of $\gamma$ to $\G$.  The map $\phi$ restricted to $\G^\circ$ is an
immersion.
\end{lem}

\begin{proof}
Let $\G$ be the linear gallery determined by a geodesic $\gamma$ and
let $\phi:\G \rightarrow K$ be canonical map.  The map $\phi$ is
clearly an embedding when restricted to each closed cell $\alpha_i$ in
$\G$.  For a point $x$ in the interior of a top cell of $\G$, there is
a neighborhood of $x$ which lies completely in $\alpha_i^\circ$.  Thus
$\phi$ is an embedding in a neighborhood of $x$.  If $x$ is a point in
the interior of a bottom cell $\alpha_i$ (which is not an endcell of
$\G$), then there is a neighborhood of $x$ which lies in
$(\alpha_{i-1} \cup_{\alpha_i} \alpha_{i+1}) \cap \G^\circ$.  The
restrictions on $\gamma$ now ensure that $\phi$ is an embedding when
restricted to this subcomplex. Finally, when $\alpha_i$ is a bottom
cell which is an endcell, there is neighborhood of $x$ which lies in
either $\alpha_{i+1} \cap \G^\circ$ or $\alpha_{i-1} \cap \G^\circ$.
\end{proof}

\begin{lem}\label{lem:retraction}
Let $K$ be an $M_\kappa$-complex, let $\phi:\G\rightarrow K$ be the
linear $[$or circular$]$ gallery determined by a local $[$closed$]$
geodesic $\gamma:[0,\ell] \rightarrow K$, and let $\gamma':[0,\ell]
\rightarrow \G$ be the lift of $\gamma$ to $\G$.  There exists a
deformation retraction from $\G^\circ$ to the image of $\gamma'$.
\end{lem}

\begin{proof}
To prove there is a deformation retraction from $\G$ to $\gamma'$ we
first note that for each bottom cell of $\G$, we can use a radial
retraction from $\alpha_i^\circ$ to the unique point where $\gamma'$
intersects $\alpha_i^\circ$ (Lemma~\ref{lem:bottom-cells}).  Next for
each top cell $\alpha_i$ there are two points which determine the
intersection of $\gamma'$ with $\alpha_i$.  These two points are
either the two endpoints of $\gamma'$, an endpoint and a point in the
interior of a boundary cell of $\alpha_i$, or two points in the
interiors of two incompatible boundary cells of $\alpha_i$.  In all
three cases it is easy to see that there is a retraction of $\alpha_i
\cap \G^\circ$ which is compatible with the deformations already
defined on the interiors of the boundary cell(s) of $\alpha_i$.
\end{proof}

\begin{defn}[Equivalent galleries]
Let $K$ be an $M_\kappa$-complex and let $\G$ and $\G'$ be two of its
galleries, both linear or both circular.  The galleries $\G$ and $\G'$
will be considered \emph{equivalent} if there is an isometry between
them which preserves their cell structures.  They will be considered
\emph{identical} if they are equivalent and they are immersed in $K$
in the same way. More specifically, if $\G$ and $\G'$ are immersed in
$K$ by maps $\phi$ and $\phi'$, then there is a map
$\rho:\G\rightarrow \G'$ showing they are equivalent with $\phi =
\phi' \circ \rho$.
\end{defn}

\begin{lem}\label{lem:finite}
Let $K$ be a $M_\kappa$-complex with $\shapes(K)$ finite, let $\ell$
be a fixed real number.  If a number $N$ satisfying the conclusion of
Theorem~\ref{thm:quasi-geo} can be constructed from $\shapes(K)$ and
$\ell$, then there exists a finite, constructible list of linear
galleries such that every geodesic in $K$ of length at most $\ell$
determines a gallery equivalent to a linear gallery in this list.
There is a similar list of circular galleries such that every closed
geodesic in $K$ of length at most $\ell$ determines a circular gallery
equivalent to a gallery in this list.
\end{lem}

\begin{proof}
We prove the linear case; the circular case is analogous.  Let
$\Gamma$ be the (infinite) collection of all geodesics in $K$ of
length at most $\ell$.  By Theorem~\ref{thm:quasi-geo} there is an $N$
such that each $\gamma$ in $\Gamma$ passes through at most $N$ open
cells of $K$.  As a result each $\gamma$ in $\Gamma$ determines a
gallery formed by gluing together at most $N$ closed cells, each of
which is isometric to a closed cell in $\shapes(K)$.  Because
$\shapes(K)$ is finite, there are only a finite number of complexes
which can be formed in this way, and by Remark~\ref{rem:recognize} we
can determine which of these are linear galleries.  The linear gallery
determined by $\gamma$ must be equivalent to one of them.
\end{proof}

%%%%%%%%%%%%%%%%%%%%%%%%%%%%%%%%%%%%%%%%%%%%%%%%%%%%%%%%%%%%
\section{Detecting piecewise geodesics}\label{sec:piecewise}
%%%%%%%%%%%%%%%%%%%%%%%%%%%%%%%%%%%%%%%%%%%%%%%%%%%%%%%%%%%%
In this section we discuss in detail how to search for piecewise
geodesics of bounded length in a linear or circular gallery using
boolean combinations of polynomial equations and inequalities.

\begin{thm}\label{thm:lin-piecewise}
If $K$ is finite $n$-dimensional $\sph$-complex and $\G$ is a linear
gallery in $K$, then we can construct a system of polynomial equations
and inequalities which has a solution if and only if the gallery $\G$
contains a piecewise geodesic of length less than $\pi$.  As a
consequence, there exists an algorithm to test whether $\G$ contains
such a piecewise geodesic.
\end{thm}

\begin{proof}
Let $\gamma:[0,\ell] \rightarrow K$ be a piecewise geodesic which is
contained in $\G$, and let $\gamma':[0,\ell]\rightarrow \G$ be its
lift to $\G$.  For each bottom cell $\alpha$ in $\G$, let $x_\alpha$
be the smallest value in the preimage of $\alpha^\circ$ under
$\gamma'$.  Such a value exists by Lemma~\ref{lem:subintervals}.
Clearly the ordering of the $x_\alpha$ is consistent with the linear
ordering of the bottom cells in $\G$.  Next, define a list $L =
\{0=x_0 < x_1 < \cdots < x_k=\ell\} \subset [0,\ell]$ which consists
endpoints of the domain together with the collection of $x_\alpha$
defined above.  Let $\alpha_1^\circ$ and $\alpha_k^\circ$ denote the
open cells of $\G$ containing $\gamma'(0)$ and $\gamma'(\ell)$ and let
$\alpha_i^\circ$, $1 < i < k$ denote the open bottom cell of $\G$
containing $\gamma'(x_i)$.

If there is an $i$ such that $\gamma'$ does not send this subinterval
$[x_i,x_{i+1}]$ to the unique geodesic connecting $\gamma'(x_i)$ and
$\gamma'(x_{i+1})$ in the unique top cell containing both points, then
we can define a new path $\gamma''$ with this property which
determines the same list $L$, but with a strictly shorter length.
Thus, if $\G$ contains a path of length less than $\pi$, it contains a
path of this type.  We may therefore assume that $\gamma'$ has this
shortest geodesic property for each $i$, and as a result that
$\gamma'$ intersects each bottom cell of $\G$ in a unique point.

We now establish a system of equations and inequalities which will
test whether the sum of the distances from $x_i$ to $x_{i+1}$,
$i=0,\ldots,k-1$ is less than $\pi$ or not.  In one sense all we need
to do is to add up the lengths of the unique individual arcs
connecting $\gamma'(x_i)$ to $\gamma'(x_{i+1})$, but we wish to do
this with polynomial equations.  Let $V$ denote the set of $0$-cells
in $\G$ and let $E$ denote the set of $1$-cells in $\G$.  We will
introduce the necessary variables and equations in three stages.  The
spaces and maps used are summarized in Figure~\ref{fig:lin-gallery}.

\begin{figure}
$\begin{array}{ccccccc}
& & [0,\ell] & \rightarrow & \sph^1 & \subset & \R^2\\
& \swarrow & \downarrow &&&&\\
K & \leftarrow & \G & \rightarrow & \sph^n & \subset & \R^{n+1}
\end{array}$
\caption{Spaces and maps used in Theorem~\ref{thm:lin-piecewise}
\label{fig:lin-gallery}}
\end{figure}

{\bf Step 1:} For every $v \in V$ we introduce a vector $\vec{u}_v =
(u_{v,1},\ldots,u_{v,n+1})$ in $\R^{n+1}$ and we add the equation
$\vec{u}_v \cdot \vec{u}_v = 1$ to our system.  In other words, we
introduce $n$ new variables $u_{v,i}$, $i=1,\ldots,n+1$ and the
polynomial equation $u_{v,1}^2 + u_{v,2}^2 +\cdots+ u_{v,n+1}^2= 1$.
The vector notation is simply a convenient shorthand.  Next, for each
$1$-cell $e \in E$ connecting $v$ to $v'$ we add the equation
$\vec{u}_v \cdot \vec{u}_{v'} = \cos(\theta)$ to our system, where
$\theta$ is the arc length of $e$.  Notice that the solutions to this
system of equations defined so far are in one-to-one correspondence
with the possible maps from $\G$ to $\sph^n$ which restrict to an
isometry on each closed cell of $\G$.

{\bf Step 2:} For each $x_i$ in $L$, we introduce a vector $\vec{y}_i
= (y_{i,1},\ldots,y_{i,n+1})$, and add the equation $\vec{y}_i\cdot
\vec{y}_i = 1$ which stipulates that $\vec{y}_i$ has length $1$.  In
addition, we add equations and inequalites which stipulate that
$\vec{y}_i$ is a positive linear combination of the vectors
$\vec{u}_v$ corresponding to the $0$-cells $v$ of $\alpha_i$.  When
$\alpha_i$ is $0$-dimensional the positive linear combination will
reduce to a set of equations of the form $\vec{y}_i = \vec{u}_v$.

Notice that we are treating the position of $\gamma(x_i)$ in
$\alpha_i^\circ$ as unspecified.  For example, if the start point
$\gamma(x)$ lies in a $1$-cell $e$ in $\G$ with $0$-cells $v$ and
$v'$, then the equations added so far only state that the image of the
endpoint $x$, $\vec{y}_1$, is a positive linear combination of the
vectors $\vec{u}_v$ and $\vec{u}_{v'}$.  They do not state what
coefficients of that positive linear combination are.  The variables
and equations in steps~$1$ and~$2$ ensure that encoded in each
solution to these system is the description of a distinct piecewise
geodesic path between the endcells of $\G$.

{\bf Step 3:} To check whether the total length of $\gamma'$ is at
least $\pi$, we use a $2$-dimensional model space.  For each point
$x_i$ in $L$ we introduce a new vector $\vec{z}_i = (z_{i,1},z_{i,2})$
in $\R^2$ starting with $\vec{z}_0 = (0,1)$.  Then we add the
following equations:
\[ \begin{array}{rcll}
\vec{z}_i \cdot \vec{z}_i &=& 1,&\textrm{ for }i=0,1,\ldots,k\\
\vec{z}_{i-1} \cdot \vec{z}_i &=& 
\vec{y}_{i-1} \cdot \vec{y}_i,&\textrm{ for }i=1,\ldots,k\\
\det(\vec{z}_{i-1},\vec{z}_i) &>& 0,&\textrm{ for }i=1,\ldots,k\\
\end{array}
\]

\noindent These equations stipulate that the length of $\vec{z}_i$ is
$1$, that the angle between $\vec{z}_i$ and $\vec{z}_{i+1}$ equals the
angle between $\vec{y}_i$ and $\vec{y}_{i+1}$, and that
$(\vec{z}_{i},\vec{z}_{i+1})$ is a positively oriented frame for
$\R^2$.  The third set of equations has this interpretation because
the cells of $\G$ are proper and thus all distances in these cells are
less than $\pi$.  Taken together, the second and third sets of
equations ensure that $\vec{z}_i$ is obtained from $\vec{z}_{i-1}$ by
a counterclockwise rotation through an angle equal to the distance
from $\gamma(x_{i-1})$ to $\gamma(x_i)$ in $\G$.  In other words the
points $\vec{z}_0$, $\vec{z}_1$, etc. are marching around the unit
circle in $\R^2$ in a counterclockwise direction starting on the
positive $x$-axis.

{\bf Step 4:} Finally, since each angle is less than $\pi$, the length
of the geodesic will be less than $\pi$ if and only if each
$\vec{z}_i$, $i>0$, has a positive second coordinate.  Thus we add the
equations $z_{i,2}>0$, for $i=1,\ldots,k$.

Combining all of the variables and equations introduced in each of
these steps, we get a system of equations and inequalities which has a
solution if and only if $\G$ contains a piecewise geodesic of length
less than $\pi$.  The second assertion follows immediately from
Tarski's theorem.
\end{proof}

\begin{figure}
\begin{tabular}{cc}
\begin{tabular}{c}\includegraphics{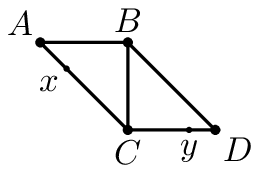}\end{tabular}
\begin{tabular}{c}\includegraphics{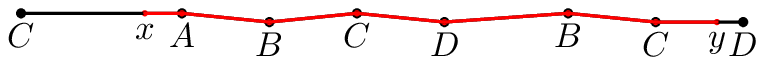}\end{tabular}
\end{tabular}
\caption{The $1$-complex and gallery described in
Example~\ref{exmp:1-gallery}
\label{fig:1-gallery}}
\end{figure}

\begin{figure}
\includegraphics{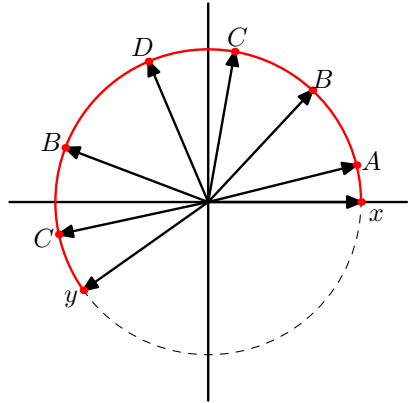}
\caption{Model space determined by $\{x,A,B,C,D,B,C,y\}$.
\label{fig:1-model}}
\end{figure}

\begin{exmp}\label{exmp:1-gallery}
To illustrate the procedure described in
Theorem~\ref{thm:lin-piecewise} consider the $1$-dimensional
$\sph$-complex on the left in Figure~\ref{fig:1-gallery}.  The gallery
determined by the sequence of points $\{x,A,B,C,D,B,C,y\}$ is shown on
the right and the model space for this sequence is shown in
Figure~\ref{fig:1-model}.  The vectors in the model space are the
vectors $\vec{z}_i$.  In this example $\vec{z}_6$ and $\vec{z}_7$ have
negative second coordinates, the total length is at least $\pi$, and
thus this is not a solution to the system of equations and
inequalities.
\end{exmp}

There is a similar theorem for circular galleries.

\begin{thm}\label{thm:cir-piecewise}
If $K$ is a finite $n$-dimensional $\sph$-complex and $\G$ is a
circular gallery in $K$, then we can construct a boolean combination
of polynomial equations and inequalities which has a solution if and
only if the gallery $\G$ contains a closed piecewise geodesic of
length less than $2\pi$.  As a consequence, there exists an algorithm
to test whether $\G$ contains such a closed piecewise geodesic.
\end{thm}

\begin{proof}
The proof is nearly identical to the one above, so only the
differences will be noted.  We begin noting that the circular gallery
$\G$ was constructed by first creating a linear gallery $\G'$ and then
identifying the two endcells of $\G'$ by an isometry.  We will need
the unidentified linear gallery $\G'$ to create some of our equations.
By picking a point in top cell as a basepoint of the loop, we may
assume that the endcells of $\G'$ are top cells.

The notation is as before, but slightly modified since $\G$ is
circular.  There is a loop $\gamma':C_\ell\rightarrow \G$ and a path
$\gamma'':[0,\ell]\rightarrow \G'$.  There is a list $L =\{0 < x_1 <
x_2 < \cdots < x_k < \ell\} \subset [0,\ell)$ where the point
$\gamma'(x_i)$ is contained in the interior of a bottom cell of $\G'$
denoted $\alpha_i$.  Moreover, $\gamma'$ is a piecewise geodesic whose
image is determined by the points $\gamma'(x_i)$.  Let $V$ denote the
set of $0$-cells in $\G'$, let $E$ denote the set of $1$-cells in
$\G'$. 

{\bf Step 1:} Same as above with $\G$ replaced by $\G'$.  The
solutions to this system of equations are in one-to-one correspondence
with the possible maps from $\G'$ to $\sph^n$ which restrict to an
isometry on each closed cell of $\G'$.

{\bf Step 2:} Same as above with the addition that since $x_1$ and
$x_k$ are two preimages of the same point in $\G$, we add equations
that state that each coefficient variable used to describe $\vec{y}_0$
as a positive linear combination of the vectors $\vec{u}_v$
corresponding to the $0$-cells in $\alpha_1$, is equal to the
corresponding coefficient variable used to describe $\vec{y}_k$ as a
positive linear combination of the vectors $\vec{u}_v$ corresponding
to the $0$-cells in $\alpha_k$.  The correspondence of the $0$-cells
is determined by the isometry used to obtain $\G$ from $\G'$.  The
variables and equations in steps~$1$ and~$2$ ensure that encoded in
any solution to the system so far is the description of a piecewise
geodesic loop which travels monotonically around $\G$.

{\bf Step 3:} Same as above.

{\bf Step 4:} When we were checking length less than $\pi$ this only
involved checking whether all the second coordinates were positive.
This time the length calculation should determine whether the length
is less than $2\pi$ which is more involved.  For each $i=
2,\ldots,k-1$ add new variables $p_i$ and $q_i$ and equations of the
form
\[ \begin{array}{c}
p_i\vec{z}_i +q_i \vec{z}_{i+1} = \vec{z}_1\\
p_i < 0 \textrm{ or } q_i <0 
\end{array}
\]

This is just a polynomial way of saying that $\vec{z}_1$ is not a
nonnegative linear combination of $\vec{z}_i$ and $\vec{z}_{i+1}$
(which ensures that the total length is less than $2\pi$).  To see
this note that if $z_{i+1}$ ever comes back around the full circle in
this model space, then $z_1$ will be a nonnegative linear combination
of $z_i$ and $z_{i+1}$, $i>1$.

Combining all of the variables and equations introduced in each of
these steps, we get a boolean combination of polynomial equations and
inequalities which has a solution if and only if $\G$ contains a
piecewise geodesic loop of length less than $2\pi$.  The second
assertion follows immediately from Tarski's theorem.
\end{proof}

\begin{rem}[Coefficients]\label{rem:coefficients}
The only equations in these systems with potentially non-rational
coefficients are the equations of the form $\vec{u}_v \cdot
\vec{u}_{v'} = \cos(\theta)$.  Such an equation can be rewritten as
the equation $(\vec{u}_v \cdot \vec{u}_{v'})^2 = \cos^2(\theta)$, and
an inequality stipulating that $\vec{u}_v \cdot \vec{u}_{v'}$ has the
appropriate sign.  Thus, if $\cos(\theta)$ is at worst a square root
of a rational number, using this alternative ensures the systems use
only rational coefficients.  Such considerations can play a
significant role when implementing algorithms, and we will explore the
issue of coefficients in greater detail in
Section~\ref{sec:coefficients}.
\end{rem}

%%%%%%%%%%%%%%%%%%%%%%%%%%%%%%%%%%%%%%%%%%%%%%%%%%%
\section{Detecting local geodesics}\label{sec:geos}
%%%%%%%%%%%%%%%%%%%%%%%%%%%%%%%%%%%%%%%%%%%%%%%%%%%
In this section we discuss how to search for local geodesics of
bounded length in a linear or circular gallery using boolean
combinations of polynomial equations and inequalities.  The
requirement that the paths be local geodesics instead of merely
piecewise geodesics is the main complication.  In the process we show
how to bound the length of a linear gallery which can contain a short
geodesic.  The results in this section will be proved by induction on
$n$.  We begin with dimension~$n=1$.

\begin{thm}\label{thm:lin-local}
If $K$ is a finite $n$-dimensional $\sph$-complex and $\G$ is a linear
gallery in $K$, then there exists a boolean combination of polynomial
equations and inequalities which has a solution if and only if $\G$
contains a local geodesic of length less than $\pi$ whose image in $K$
is also a geodesic.  As a consequence there exists an algorithm which
determines whether or not $\G$ contains such a geodesic.
\end{thm}

\begin{proof}
When $n=1$, the system created is the one created in the proof of
Theorem~\ref{thm:lin-piecewise}.  Because $K$ and hence $\G$ have
dimension~$1$, the piecewise geodesics encoded in a solution to this
system will be precisely the local geodesics contained in $\G$, and
they will also remain local geodesics when immersed into $K$.

For $n>1$, we may assume that results in this section have already
been proved in lower dimensions.  Let $P$ denote the system created in
the proof of Theorem~\ref{thm:lin-piecewise}. In essence, we will
start with $P$ and then add restrictions to require that any piecewise
geodesics encoded in a solution to $P$ will be locally geodesic at the
transitions between the geodesic pieces.  We will use the notations
established in the proof of Theorem~\ref{thm:lin-piecewise} without
further comment.

Let $\alpha_i$ be a bottom cell of $\G$.  Given any two points $a$ and
$b$ in $\alpha_i^\circ$, the links $\link(\phi(a),K)$ and
$\link(\phi(b),K)$ are canonically isometric.  Let $L_i$ be this
well-defined $\sph$-complex and note that its dimension is strictly
less than $n$.  In general the cells in $L_i$ will not be proper so we
assume that it has been suitably subdivided.  Notice that the three
points $\gamma(x_{i-1})$, $\gamma(x_{i})$ and $\gamma(x_{i+1})$
determine two points in $L_i$.  The complex $L_i$ will have a finite
number of galleries which could possibly contain a geodesic of length
less than $\pi$ by Lemma~\ref{lem:pi}.  Let $G_i$ be the list of such
galleries in $L_i$.  For each gallery in $G_i$ we write down the
variables and equations necessary to the encode a generic path
determined by crossing points into the gallery (this would be steps
$1$ and $2$ in Theorem~\ref{thm:lin-piecewise}).  Do this for each
gallery in $G_i$ and repeat this procedure for each of the links
$L_i$.  Each gallery in each link will use distinct variables.

For each $i$ and for each gallery $\script{H}_j$ in $G_i$, we create a
system $V_{i,j}$ which is a boolean combination of polynomial
equations and inequalities.  First we create the equations and
inequalities needed to say that the three points $\gamma(x_{i-1})$,
$\gamma(x_{i})$ and $\gamma(x_{i+1})$ determine two points in $L_i$
which lie in the endcells of $\script{H}_j$.  Then we add the boolean
combination of polynomial equations and inequalities (which exist by
the lower-dimensional versions of this theorem) which have a solution
if and only if there is a local geodesic in $\script{H}_j$ connecting
these exact points in its endcells.

Intuitively, the system $V_{i,j}$ has a solution if and only if the
the piecewise geodesic in $\G$ which is encoded in a solution of the
system $P$, there is a path less than $\pi$ contained in the gallery
$\script{H}_j$ in the link of the point $\gamma(x_i)$ pushed into $K$.
More informally, we might say that the solutions to $V_{i,j}$ encode
the piecewise geodesics in $\G$ which have a ``kink'' at
$\phi(\gamma(x_i))$ as witnessed by a short path through gallery
$\script{H}_j$.

The final system is one which starts with $P$ and subtracts off the
solutions to each of the $V_{i,j}$.  Clearly, this is a boolean
combination of these solution sets.  This final system will have a
solution if and only if $\G$ contains a piecewise geodesic of length
less than $\pi$ which does not fail to be locally geodesic in any of
the finite number of possible ways which are available to it.  In
other words, if and only if $\G$ contains a local geodesic of length
less than $\pi$.  Finally, the algorithmic assertion follows
immediately from Tarski's Theorem.
\end{proof}

There is a similar theorem for circular galleries.

\begin{thm}\label{thm:cir-gallery}
If $K$ is a finite $n$-dimensional $\sph$-complex and $\G$ is a
circular gallery in $K$, then there exists a boolean combination of
polynomial equations and inequalities which has a solution if and only
if $\G$ contains a local geodesic of length less than $2\pi$ whose
image in $K$ is also a local geodesic.  As a consequence there exists
an algorithm which determines whether or not $\G$ contains such a
geodesic.
\end{thm}

\begin{proof}
The proof is identical to the one above, but starting with the system
created in the proof of Theorem~\ref{thm:cir-piecewise}.  Everything
else is unchanged.
\end{proof}

Next, we show how to construct the number $N$ in dimension $n$.

\begin{lem}\label{lem:N}
If $K$ be a $n$-dimensional $\sph$-complex with finitely many shapes
and $\ell>0$ is a real number, then there is an algorithm which
computes an integer $N$ such that every local geodesic contained in a
linear geodesic $\G$ of size $N$ in $K$ has length at least $\ell$.
\end{lem}

\begin{proof}
Note that Theorem~\ref{thm:lin-local} shows that it is possible to
determine whether any specific linear gallery $\G$ contains a local
geodesic of length less than $\pi$.  To find the smallest value of $N$
which will work for $\ell = \pi$ we simply test a representative of
each equivalence class of galleries of size~$1$, then a representative
of each equivalence class of galleries of size~$2$, etc.  That there
is a finite constructible list of such galleries is guaranteed by
Lemma~\ref{lem:finite}. As soon as we find an $N$ such that all
galleries of size~$N$ do not contain a local geodesic of less than
length~$\pi$, we are essentially done since by Remark~\ref{rem:N} it
is sufficient to find $N$ for a single number $\ell$.  That such an
$N$ must exist follows from Theorem~\ref{thm:quasi-geo} and the fact
that local geodesics contained in linear galleries of size $N$
themselves have size $N$.  Note that we are using the fact that
geodesics ``contained in'' $\G$ remain in the interior of the gallery
(Remark~\ref{rem:geos-in-galleries}).  A geodesic between the endcells
of $\G$ which is allowed to contain points from the boundary of $\G$
may very well have a smaller size.
\end{proof}

Now that we have shown that the $N$ in Theorem~\ref{thm:quasi-geo} is
constructible, we can state the following corollaries.

\begin{cor}\label{cor:locally-finite}
Let $K$ be a finite $n$-dimensional $\sph$-complex and let $\ell$ be a
fixed real number.  Given open cells $\sigma^\circ$ and $\tau^\circ$
in $K$ there is a finite, constructible list of linear galleries from
$\sigma$ to $\tau$, based only on $\shapes(K)$ and $\ell$, such that
every geodesic from $\sigma$ to $\tau$ of length at most $\ell$
determines a gallery in this list.
\end{cor}

\begin{proof}
Let $N$ be a size which guarantees a length greater than $\ell$ and
then construct all possible linear galleries with size at most $N$ up
to equivalence (Lemma~\ref{lem:finite}), and then construct all
possible maps from $\G$ to $K$ starting at $\sigma$.  This is possible
since $N$ exists (Theorem~\ref{thm:quasi-geo}), it is constructible
(Lemma~\ref{lem:N}), and $K$ is locally-finite.
\end{proof}

\begin{cor}\label{cor:cir-finite}
Let $K$ be a finite $n$-dimensional $\sph$-complex and let $\ell$ be a
fixed real number.  There is a finite, constructible list of circular
galleries such that every closed geodesic of length at most $\ell$
determines a circular gallery in this list.
\end{cor}

\begin{proof}
Let $N$ be a size which guarantees a length greater than $\ell$ and
then construct all possible circular galleries $\G$ with size at most
$N$ up to equivalence (Lemma~\ref{lem:finite}), and then construct all
possible maps from $\G$ to $K$.  This is possible since $N$ exists
(Theorem~\ref{thm:quasi-geo}), it is constructible
(Lemma~\ref{lem:N}), and $K$ is locally finite.
\end{proof}

Finally, we prove a lemma to complete the induction.

\begin{lem}\label{lem:pi}
Let $K$ be a finite $n$-dimensional $\sph$-complex.  Given open cells
$\sigma^\circ$ and $\tau^\circ$ in $K$, there exists a boolean
combination of polynomial equations and inequalities which has a
solution if and only if there is a geodesic in $K$ from some point in
$\sigma^\circ$ to some point in $\tau^\circ$ of length less than
$\pi$.  A similar statement holds when the points in $\sigma^\circ$
and $\tau^\circ$ are specified.
\end{lem}

\begin{proof}
First note that any such geodesic determines a linear gallery. By
Corollary~\ref{cor:locally-finite} we can construct a finite list $L$
which contains all linear galleries $\G\rightarrow K$ from $\sigma$ to
$\tau$ which could possibly be determined by a geodesic of length less
than $\pi$.  For each gallery in $L$, use Theorem~\ref{thm:lin-local}
to create a system of polynomial equations which has a solution if and
only if this particular gallery does not contain a piecewise geodesic
of length less than $\pi$.  If each of these systems are created using
distinct variables, then their union will be a finite system which has
a solution if and only if there does not exist a path in $K$ from
$\sigma^\circ$ to $\tau^\circ$ of length less than $\pi$.
\end{proof}

%%%%%%%%%%%%%%%%%%%%%%%%%%%%%%%%%%%%%%%%%%%%%%
\section{Coefficients}\label{sec:coefficients}
%%%%%%%%%%%%%%%%%%%%%%%%%%%%%%%%%%%%%%%%%%%%%%
The results in this section are a slight digression.  Let $K$ be a
finite $n$-dimensional $\euc$-complex.  In this section we show that
the type of numbers which show up as coefficients of the polynomials
used to test whether $K$ is locally $\cat(0)$ are no worse than the
types of numbers which show up as lengths of $1$-cells in $K$.  A
precise version of this claim is contained in
Theorem~\ref{thm:coefficients} where it is presented in terms of
$F$-simplices.

\begin{defn}[$F$-simplices]\label{def:F-simplices}
Let $F$ be a subfield of $\R$ and let $\sigma$ be $n$-dimensional
Euclidean simplex.  If the square of the length of each $1$-cell in
$\sigma$ lies in $F$, then $\sigma$ is called an \emph{$F$-simplex}.
\end{defn}

\begin{lem}\label{lem:edge2dotprods}
Let $F$ be a subfield of $\R$, let $\sigma \subset \R^n$ be an
$k$-dimensional Euclidean simplex with vertices $v_i$, $i=0,\ldots,k$,
and let $\vec{v}_i$ be the vector from $v_0$ to $v_i$.  The following
conditions are equivalent:
\begin{enumerate}
\item $\sigma$ is an $F$-simplex,
\item $\vec{v}_i \cdot \vec{v}_j \in F$ for all $i,j$,
\end{enumerate}
\end{lem}

\begin{proof}
($2 \Rightarrow 1$) Since $||\vec{v}_i||^2 = \vec{v}_i \cdot
\vec{v}_i$ we conclude that the square of the length of the edge
connecting $v_0$ to $v_i$ lies in $F$, and since
\[||\vec{v}_j - \vec{v}_i||^2 = (\vec{v}_j - \vec{v}_i)
\cdot (\vec{v}_j - \vec{v}_i) = ||\vec{v}_i||^2 + ||\vec{v}_j||^2 -2
\vec{v}_i \cdot \vec{v}_j\]
\noindent 
we conclude that square of the length of the edge connecting $v_i$ to
$v_j$, $i,j >0$, lies in $F$ as well.  ($1 \Rightarrow 2$) Since the
square of the length of the edge connecting $v_0$ to $v_i$ is
$||\vec{v}_i||^2 = \vec{v}_i \cdot \vec{v}_i$, we conclude that
$\vec{v}_i \cdot \vec{v}_j$ lies in $F$ when $i=j$.  Similarly, the
square of the length of the edge connecting $v_i$ to $v_j$ is
$||\vec{v}_j - \vec{v}_i||^2 = ||\vec{v}_i||^2 + ||\vec{v}_j||^2 -2
\vec{v}_i \cdot \vec{v}_j$.  Since this result lies in $F$, $-2$ is in
$F$, and the first two of its terms lies in $F$ by assumption, we
conclude that $\vec{v}_i \cdot \vec{v}_j \in F$ for $i \neq j$ as
well.
\end{proof}

%%%%%%%%

\begin{lem}\label{lem:coordinates}
Let $F$ be a subfield of $\R$ and let $\sigma \subseteq \R^n$ be an
$k$-dimensional $F$-simplex with one vertex, $v_0$, at the origin.
There is an orthogonal basis $U$ of $\R^n$ such that the squared
length of each basis vector in $U$ lies in $F$ and all of the
coordinates of all of the vertices of $\sigma$, relative to $U$, lies
in $F$.
\end{lem}

\begin{proof}
By replacing $\R^n$ with the subspace spanned by the vectors
$\vec{v}_i$ if necessary, we may assume that $k=n$ and that collection
of vectors $V = \{\vec{v}_i\}_{i=1}^n$ is a basis for $\R^n$.  Next
apply the Gramm-Schmidt process to $V$ without normalizing the final
result.  Specifically, let $\vec{u}_1 = \vec{v}_1$, and for $i >1$
inductively define
\[\vec{u}_i = \vec{v}_i - \sum_{j=1}^{i-1} 
\frac{\vec{v}_i\cdot\vec{u}_j}{\vec{u}_j\cdot\vec{u}_j} \vec{u}_j.
\]
\noindent An easy induction shows that $\vec{u}_i$ is an $F$-linear
combination of the vectors $\vec{v}_j$, $j \leq i$ so that the
fraction in front of the vector $\vec{u}_j$ in the above equation is
an element of $F$.  In particular, the change of basis matrix
connecting $U$ and $V$ is an upper triangular matrix over $F$ with
$1$'s down the diagonal.  As a result, the same is true about the
change of basis matrix going in the other direction, thus showing that
each $\vec{v}_i$ is an $F$-linear combination of the vectors
$\vec{u}_j$ with $j \leq i$.
\end{proof}

\begin{lem}\label{lem:dihedral}
If $F$ be a subfield of $\R$ and $\sigma$ is a $F$-simplex, then the
square of each trigonometric function of each dihedral angle in
$\sigma$ lies in $F$.
\end{lem}

\begin{proof}
First note that it is sufficient to prove this for the square of the
cosine of such an angle. Let $V = \{v_0,\ldots,v_k\}$ be the set of
vertices of $\sigma$.  A dihedral angle between two faces of $\sigma$
is determined by a subset $S \subset V$ and distinct $s,t \in V
\setminus S$.  In particular, there is a dihedral angle $\theta$
between the face determined by $S \cup \{s\}$ and the face determined
by $S \cup \{t\}$.  By renumbering if necessary we may assume that $S
= \{v_0,\ldots,v_{i-1}\}$, $s = v_i$, and $t = v_{i+1}$.  Using the
orthogonal basis $U$ derived in Lemma~\ref{lem:coordinates}, it is
clear that $\cos^2(\theta)$ lies in $F$.  Specifically, $\theta$ is
angle between the projections of $vec{v}_i$ and $\vec{v}_{i+1}$ into
the subspace perpendicular to the subspace determined by $S$.  In the
basis $U$, this projection is accomplished by simply ``zeroing out''
the first $i-1$ coordinates.  As a result these projections are still
$F$-linear combinations of the basis vectors.  That $\cos^2(\theta)$
lies in $F$ now follows immediately.
\end{proof}

One special case which deserves to be highlighted is where the field
$F = \Q$.  An angle $\theta$ for which $\cos^2(\theta)$ is rational is
called \emph{geodetic}.  See \cite{CoRaSa} for a computational
description of the space of geodetic angles.

\begin{cor}\label{cor:geodetic}
If $\sigma$ be a Euclidean $\Q$-simplex, then every dihedral angle in
$\sigma$ is geodetic.
\end{cor}

Notice that $\Q$-simplices are quite common.  For example, every
simplex in a simplicial decomposition of a rational polytope will be a
$\Q$-simplex.

\begin{thm}[Coefficients]\label{thm:coefficients}
Let $K$ be a finite simplicial $\euc$-complex and let $F$ be a
subfield of $\R$.  If every simplex in $K$ is an $F$-simplex, then the
systems of polynomial equations and inequalities used to test whether
$K$ is locally $\cat(0)$ have coefficients which belong to $F$.
\end{thm}

\begin{proof}
As we mentioned in Remark~\ref{rem:coefficients}, the only
coefficients used which are not integers are equal to $\cos^2(\theta)$
for some dihedral angle $\theta$ in some simplex of $K$.  But by
Lemma~\ref{lem:dihedral}, this number lies in $F$.
\end{proof}

Similar theorems could clearly be derived for other curvatures, but
they will not be pursued here.

%%%%%%%%%%%%%%%%%%%%%%%%%%%%%%%%%%%%%%%%%%
\section{The Main Theorem}\label{sec:main}
%%%%%%%%%%%%%%%%%%%%%%%%%%%%%%%%%%%%%%%%%%
At this point the proof of Theorem~\ref{thm:main} is nearly immediate.

\begin{thm}\label{thm:cat1}
There exists an algorithm which determines whether or not a finite
$\sph$-complex contains a closed geodesic of length less than $2\pi$,
and as a consequence there exists an algorithm which determines
whether or not a finite $\sph$-complex is globally $\cat(1)$.
\end{thm}

\begin{proof}
A closed geodesic of length less than $2\pi$ in the complex would
determine a circular gallery which contains it.  By
Corollary~\ref{cor:cir-finite} there is a finite, constructible list
of circular galleries in $K$ which could possibly contain a closed
geodesic of length less than $2\pi$, and by
Theorem~\ref{thm:cir-gallery} we can test each one for the presence of
a short closed geodesic which survives in $K$.  This proves the first
assertion.

To prove the second let $K$ be the finite $\sph$-complex under
consideration.  By Theorem~\ref{thm:local-to-global} and
Theorem~\ref{thm:local-via-geos} it is sufficient to show that $K$
contains no closed geodesic of length less than $2\pi$ and that for
each cell $B$ in $K$, $\link(B,K)$ contains no closed geodesic of
length less than $2\pi$.  Since $K$ and $\link(B,K)$ are finite
$\sph$-complexes for all cells $B$, the first part of the proof shows
that these are checkable conditions.
\end{proof}

\renewcommand{\thethm}{\ref{thm:main}}
\begin{thm}
There exists an algorithm which determines whether or not a finite
$M_\kappa$-complex is locally $\cat(\kappa)$.
\end{thm}
\renewcommand{\thethm}{\thesection.\arabic{thm}}
\addtocounter{thm}{-1}

\begin{proof}
For each vertex $v$ in $K$, the link $\link(v,K)$ is a finite
$\sph$-complex and by Theorem~\ref{thm:cat1} we can check whether it
is globally $\cat(1)$.  By Theorem~\ref{thm:local-via-geos} this is
sufficient.
\end{proof}

As can be seen from the proof, it is actually sufficient for the
complex $K$ to be locally-finite, so long as (1) only a finite list of
finite $\sph$-complexes occur as links of $0$-cells up to isometry,
and (2) there is a constructive procedure for enumerating this finite
list.

Computational real algebraic geometry is able to show that testing
curvature conditions is decidable and algorithmic, but the real
algebraic sets described here quickly reach the limit of what is
realistically computable using current techniques.  For the working
geometric group theorist who would like software which determines the
local curvature of a finite $M_\kappa$-polyhedral complex, there are
essentially two main options:

\begin{itemize}
\item use sophisticated techniques from computational real algebraic
geometry to speed up the computations, or
\item develop alternative, simpler algorithms which work directly with
the geometry.
\end{itemize}

Both strategies are currently under investigation.  The latter
strategy, in particular, has already produced results in the form of a
direct geometric algorithm to determine the local curvature of a
$3$-dimensional $M_\kappa$-complex using only elementary
$3$-dimensional geometry.  This result can be found in
\cite{ElMc-three}.  At this point, the following question remains
open:

\begin{prob}
Is there a direct \emph{geometric} algorithm which determines the
local curvature of a finite $M_\kappa$-complex in dimensions greater
than $3$?
\end{prob}

%%%%%%%%%%%%%%%%%%%%%%%%%%%

%%%%%%%%%%%%%%

\begin{thebibliography}{10}

\bibitem{Ba}
W.~Ballmann.
\newblock {\em Lectures on spaces of nonpositive curvature}.
\newblock Birkh\"auser Verlag, Basel, 1995.
\newblock With an appendix by Misha Brin.

\bibitem{BoCoRo}
J.~Bochnak, M.~Coste, and M.~Roy.
\newblock {\em Real Algebraic Geometry}.
\newblock Springer-Verlag, Berlin, 1998.

\bibitem{Br-thesis}
M.~R.~Bridson.
\newblock Geodesics and curvature in metric simplicial complexes.
\newblock In {\em Group theory from a geometrical viewpoint (Trieste, 1990)},
  pages 373--463. World Sci. Publishing, River Edge, NJ, 1991.

\bibitem{BrHa}
M.~R.~Bridson and A.~Haefliger.
\newblock {\em Metric spaces of non-positive curvature}.
\newblock Springer-Verlag, Berlin, 1999.

\bibitem{CaJo}
B.~F.~Caviness and J.~R.~Johnson, eds.
\newblock {\em Quantifier eliminations an cylindrical algebraic
decomposition}.
\newblock Springer-Verlag, Berlin, 1998.

\bibitem{CoRaSa}
J.~H. Conway, C.~Radin, and L.~Sadun.
\newblock On angles whose squared trigonometric functions are rational.
\newblock {\em Discrete Comput. Geom.}, 22(3):321--332, 1999.

\bibitem{ElMc-three}
M.~Elder and J.~McCammond.
\newblock Curvature testing in $3$-dimensional metric polyhedral complexes.
\newblock Preprint 2001.

\bibitem{Gr}
M.~Gromov.
\newblock Hyperbolic groups.
\newblock In {\em Essays in group theory}, pages 75--263. Springer, New York,
  1987.

\bibitem{Mo}
G.~Moussong.
\newblock Hyperbolic {C}oxeter groups.
\newblock Ph.D. dissertation, The Ohio State University, 1988.

\end{thebibliography}
\end{document}